\documentclass[11pt]{amsart}

\usepackage{amsmath}
\usepackage{amssymb}
\usepackage{amsthm}
\usepackage{enumerate}
\usepackage{amsbsy}
\usepackage{amsfonts}

\usepackage[bitstream-charter]{mathdesign}
\usepackage[T1]{fontenc}

\usepackage{amsthm}
\newtheorem{theorem}{Theorem}[section]
\newtheorem{lemma}[theorem]{Lemma}
\theoremstyle{definition}

\theoremstyle{remark}

\providecommand{\norm}[1]{\lVert#1\rVert}

\providecommand{\Bignorm}[1]{\Bigl\lVert#1\Bigr\rVert}

\begin{document}
\baselineskip=18pt

\title[Fractional Laplacian operator]{A parabolic Triebel-Lizorkin space estimate for the fractional Laplacian operator}

\author{Minsuk Yang}
\address{School of Mathematics, Korea Institute for Advanced Study, 85 Hoegiro Dongdaemungu, Seoul, Republic of Korea}
\email{yangm@kias.re.kr}

\begin{abstract}
In this paper we prove a parabolic Triebel-Lizorkin space estimate for the operator given by 
\[T^{\alpha}f(t,x) = \int_0^t \int_{{\mathbb R}^d} P^{\alpha}(t-s,x-y)f(s,y) dyds,\]
where the kernel is 
\[P^{\alpha}(t,x) = \int_{{\mathbb R}^d} e^{2\pi ix\cdot\xi} e^{-t|\xi|^\alpha} d\xi.\]
The operator $T^{\alpha}$ maps from $L^{p}F_{s}^{p,q}$ to $L^{p}F_{s+\alpha/p}^{p,q}$ continuously.
It has an application to a class of stochastic integro-differential equations of the type $du = -(-\Delta)^{\alpha/2} u dt + f dX_t$.
\end{abstract}

\maketitle

\section{Introduction}
\label{S1}

This paper is concerned with the regularity of solutions of the following stochastic partial differential equations 
\[\begin{cases}
du = -(-\Delta)^{\alpha/2} u dt + f dX_t &  (0,T) \times {{\mathbb R}^d} \\
u = 0 & \{0\} \times {{\mathbb R}^d},
\end{cases}\]
where $0<\alpha\le 2$ and $X_{t}$ are L\'evy processes.
The solution can be represented by the stochastic integral
\[u(t,x,\omega) = \int_0^t \int_{{\mathbb R}^d} P^{\alpha}(t-s,x-y) f(s,y,\omega) dydX_s(\omega),\]
where the kernel is given by 
\begin{equation}
\label{s1:e3}
P^{\alpha}(t,x) = \int_{{\mathbb R}^d} e^{2\pi ix\cdot\xi} e^{-t|\xi|^\alpha} d\xi.
\end{equation}
For general background of L\'evy processes and stochastic calculus, see, for example, Applebaum's monograph \cite{MR2512800}.

The concepts of solutions and $L^{p}$-theory of the stochastic partial differential equations were already established by N. V. Krylov \cite{MR1661766,MR2227229}.
He studied the regularity of solutions to the above problem when $\alpha=2$ and $X_{t}$ is a Brownian motion in \cite{MR1317805}.
The Burkholder--Davis--Gundy inequality implies that 
\begin{align*}
&{\mathbb E} \int_0^T\norm{\nabla_{x} u}_{L^p({\mathbb R}^d)}^p dt \\
&\ll {\mathbb E} \int_0^T \int_{{\mathbb R}^d} \Big(\int_0^t \Big|\nabla_{x} \int_{{\mathbb R}^d} P^{2}(t-s,x-y)f(s,y,\omega) dy\Big|^2 ds\Big)^{p/2} dxdt,
\end{align*}
where ${\mathbb E}X$ denotes the expectation $\int_{\Omega} X(\omega) d{\mathbb P}(\omega)$.
He proved that the above parabolic square function is bounded for $2\le p<\infty$ by 
\[{\mathbb E} \int_0^T \|f\|^p_{L^p({\mathbb R}^d)} ds,\]
interpolating an $L^{2}$-estimate via Plancherel's theorem and a bounded mean oscillation estimate.
Notice that the desired result depends on the deterministic estimates after applying the Burkholder--Davis--Gundy inequality.

	Recently, there are many studies about more general stochastic partial differential equations with the fractional Laplacian operator \cite{MR2283957, MR2946443, MR2869738} or about more general function spaces \cite{MR3119988, MR3034605}, for example.
For the L\'evy processes $X_{t}$ and $0<\alpha<2$, a few results are known for general Sobolev estimates.
In this case, the Kunita inequality applies under some condition on L\'evy processes that 
\begin{align*}
&{\mathbb E} \int_0^T \norm{\nabla u(s, \cdot)}_{L^p({\mathbb R}^d)}^p ds \\
&\ll {\mathbb E} \int_0^T \int_{{\mathbb R}^d} \Big|\nabla_{x} \int_{{\mathbb R}^d} P^{\alpha}(t-s,x-y)f(s,y) dy\Big|^2 ds\Big)^{p/2} dxdt \\
&\quad + {\mathbb E} \int_0^T \int_{{\mathbb R}^d} \int_0^t \Big|\nabla_{x} \int_{{\mathbb R}^d} P^{\alpha}(t-s,x-y)f(s,y) dy\Big|^p dsdxdt.
\end{align*}
We note that the desired result also can be deduced from certain deterministic estimates after applying the Kunita inequality.

The purpose of this paper is to further elucidate the main estimates of this kind of regularity theorems.
More precisely, we shall prove the following theorem.

\begin{theorem}
Let $Q=[0,T]\times{{\mathbb R}^d}$ with $0<T<\infty$.
If $0<\alpha\le2$, $2\le q\le p<\infty$, and $s\in{\mathbb R}$, then the operator $T^{\alpha}$ given by 
\begin{equation}
\label{s1:e2}
T^{\alpha}f(t,x) = \int_0^t \int_{{\mathbb R}^d} P^{\alpha}(t-s,x-y)f(s,y) dyds
\end{equation}
maps from $L^{p}F_{s}^{p,q}(Q)$ to $L^{p}F_{s+\alpha/p}^{p,q}(Q)$ and satisfies 
\begin{equation}
\label{s1:e1}
\norm{T^{\alpha}f}_{L^{p}F_{s+\alpha/p}^{p,q}(Q)} \ll \norm{f}_{L^{p}F_{s}^{p,q}(Q)}.
\end{equation}
Here, the spaces $L^{p}F_{s}^{p,q}(Q)$ is the set of measurable functions $u:Q\to{\mathbb R}$ such that 
\[\norm{u}_{L^{p}F_{s}^{p,q}(Q)} := \left(\int_{0}^{T} \norm{u(t,\cdot)}_{F_{s}^{p,q}({{\mathbb R}^d})}^{p} dt\right)^{1/p} < \infty\]
and $F_{s}^{p,q}({{\mathbb R}^d})$ is the standard Triebel-Lizorkin space.
\end{theorem}

Throughout the paper, we shall use the notation $a \ll b$, which means $|a| \le cb$ for some positive number $c$. 

Our method of proof was based on the Littlewood--Paley theory, Khinchine's inequality, and weighted convexity inequalities.
By the direct application of the standard vector-valued singular integral theory, it is not easy to obtain this regularity result.
Our proof is simple because we use the decay of the operator norms of the Littlewood--Paley pieces effectively.

\section{Preliminaries}

In this section, we set down notations and definitions.

\subsection{Triebel--Lizorkin spaces}

Given a Schwartz function $f$, we define its Fourier transform by 
\[\widehat{f}(\xi) := \int_{{\mathbb R}^d} e^{-2\pi i \xi \cdot x} f(x) dx.\]
The definition of Fourier transform extends naturally to tempered distributions.

It is a remarkable fact that several function spaces are characterized by using Littlewood--Paley theory. 
To define the Littlewood--Paley operators, we fix a radial Schwartz function $\Phi$ on ${{\mathbb R}^d}$ whose Fourier transform is nonnegative, supported in the ball $|\xi|\le2$, equal to 1 in the ball $|\xi|\le1$ and define $\widehat{\Psi}(\xi)=\widehat{\Phi}(\xi)-\widehat{\Phi}(2\xi)$.
We have the following partitions of unity
\begin{equation}
\label{s2:e1}
\widehat{\Phi}(\xi)+{\sum_{j=1}^\infty} \widehat{\Psi}(2^{-j}\xi)=1.
\end{equation}
Then the Littlewood-Paley operators $S_0$ and $\Delta_j$ for all integers $j$ are defined as 
\begin{equation}
\label{s2:e2}
S_0 f(x) := \int_{{\mathbb R}^d} e^{2\pi i x\cdot\xi} \widehat{\Phi}(\xi)\widehat{f}(\xi) d\xi
\end{equation}
and 
\begin{equation}
\label{s2:e3}
\Delta_j f(x) := \int_{{\mathbb R}^d} e^{2\pi i x\cdot\xi} \widehat{\Psi}(2^{-j}\xi)\widehat{f}(\xi) d\xi.
\end{equation}

Now, we define the Triebel-Lizorkin spaces.
Let $s\in{\mathbb R}$ and $0<p,q\le\infty$.
The Triebel-Lizorkin space $F_{p,q}^s$ is the space of all tempered distributions $f$ with 
\begin{equation}
\label{s2:e4}
\norm{f}_{F_{s}^{p,q}} := \norm{S_0f}_p + \Bignorm{\Big({\sum_{j=1}^\infty} (2^{js}|\Delta_jf|)^q\Big)^{1/q}}_p < \infty.
\end{equation}

We remark that for $1<p<\infty$ and $s\in{\mathbb R}$, two Banach spaces $F_{s}^{p,2}$ and $L_{s}^{p}$ have equivalent norms.
So, the Triebel-Lizorkin spaces are the natural generalization of the fractional Sobolev spaces.

\subsection{Khinchine's inequality}

The Rademacher functions $r_{1}(z), r_{2}(z), \dots, r_{j}(z), \dots$ are defined on the interval $[0,1]$ as follows:
\[\begin{cases}
r_{1}(z)=1 &  \text{ for } 0 \le z \le 1/2 \\
r_{1}(z)=-1 &  \text{ for } 1/2 < z <1.
\end{cases}\]
It is extended outside the unit interval by period 1.
In general 
\[r_{j}(z)=r_{1}(2^{j-1}z).\]

The Rademacher functions provide a very useful device in the study of Lebesgue norms in terms of quadratic expressions.

\begin{lemma}[Khinchine's inequality]
There exists a positive constant $C$ depending only on $p$ such that for any sequence of complex numbers $c_{j}$,
\begin{equation}
\label{s2:e5}
\frac{1}{C} \Big({\sum_{j=1}^\infty} |c_{j}|^{2}\Big)^{p/2} \le \int_{0}^{1} \Big|{\sum_{j=1}^\infty} c_{j}r_{j}(z)\Big|^{p} dz \le C \Big({\sum_{j=1}^\infty} |c_{j}|^{2}\Big)^{p/2}.
\end{equation}
\end{lemma}

This is a consequence of sub-Gaussian bounds.
It follows from the fact that the sequence of Rademacher functions are mutually independent over $[0,1]$ and take values $\pm 1$ with equal probability.
For the detailed proof of Khinchine's inequality, see the appendix of Stein's monograph \cite{MR0290095}.

\section{Proof of the theorem}

\textbf{Step 1:} 
First we note that 
\[\norm{T^{\alpha}f}_{L^p F_{s+\alpha/p}^{p,q}}^{p} \ll \int_{0}^{T} \int_{{\mathbb R}^d} |S_{0}T^{\alpha}f|^{p} dxdt + \int_{0}^{T} \int_{{\mathbb R}^d} \Big({\sum_{j=1}^\infty} |2^{j(s+\alpha/p)} \Delta_{j}T^{\alpha}f|^{q}\Big)^{p/q} dxdt.\]
Because $S_{0}T^{\alpha}f=T^{\alpha}S_{0}f$ and the kernel \eqref{s1:e3} is integrable in $x$ uniformly in $t$, we have by Young's inequality
\[\int_{0}^{T} \int_{{\mathbb R}^d} |S_{0}T^{\alpha}f|^{p} dxdt \ll \int_{0}^{T} \int_{{\mathbb R}^d} |S_{0}f|^{p} dxdt.\]
Let $r_{j}(z)$ denote the Rademacher functions.
For $2\le q$, using Khinchine's inequality \eqref{s2:e5}, we obtain that  
\begin{align*}
&\int_{0}^{T} \int_{{\mathbb R}^d} \Big({\sum_{j=1}^\infty} |2^{j(s+\alpha/p)} \Delta_{j}T^{\alpha}f(t,x)|^{q}\Big)^{p/q} dxdt \\
&\le \int_{0}^{T} \int_{{\mathbb R}^d} \Big({\sum_{j=1}^\infty} |2^{j(s+\alpha/p)} \Delta_{j}T^{\alpha}f(t,x)|^{2}\Big)^{p/2} dxdt \\
&\ll \int_{0}^{T} \int_{{\mathbb R}^d} \int_{0}^{1} \Big|{\sum_{j=1}^\infty} r_{j}(z) 2^{j(s+\alpha/p)} \Delta_{j}T^{\alpha}f(t,x)\Big|^{p} dzdxdt.
\end{align*}
Now, we observe that 
\[\Delta_jT^{\alpha}f(t,x) = \Delta_jT^{\alpha}(\Delta_{j-1}+\Delta_j+\Delta_{j+1})f(t,x).\]
This is easily verified by taking the Fourier transform with respect to the space variable.
So, we have 
\begin{equation}
\label{s3:e1}
\begin{split}
&\int_{0}^{T} \int_{{\mathbb R}^d} \int_{0}^{1} \Big|{\sum_{j=1}^\infty} r_{j}(z) 2^{j(s+\alpha/p)} \Delta_{j}T^{\alpha}f(t,x)\Big|^{p} dzdxdt \\
&\ll \int_{0}^{T} \int_{{\mathbb R}^d} \int_{0}^{1} \Big|{\sum_{j=1}^\infty} r_{j}(z) 2^{j(s+\alpha/p)} \Delta_{j}T^{\alpha}\Delta_{j}f(t,x)\Big|^{p} dzdxdt.
\end{split}
\end{equation}
Let us denote 
\[F_{j}(s,y) = 2^{js}\Delta_{j}f(s,y)\]
and
\[P_{j}^{\alpha}(t,x) = 2^{j\alpha/p} \int_{{\mathbb R}^d} e^{2\pi ix\cdot\xi} \widehat{\Psi}(2^{-j}\xi) e^{-t|\xi|^\alpha} d\xi.\]
Then we can write the operator in \eqref{s3:e1} as 
\[2^{j(s+\alpha/p)} \Delta_{j}T^{\alpha}\Delta_{j}f(t,x) = \int_{0}^{t} \int_{{\mathbb R}^d} P_{j}^{\alpha}(t-s,x-y) F_{j}(s,y) dy ds.\]
Therefore, using Jensen's inequality and the triangle inequality, we can write the last integral in \eqref{s3:e1} as 
\begin{equation}
\label{s3:e2}
\begin{split}
&\int_{0}^{T} \int_{{\mathbb R}^d} \int_{0}^{1} \Big|\int_{0}^{t} {\sum_{j=1}^\infty} r_{j}(z) \int_{{\mathbb R}^d} P_{j}^{\alpha}(t-s,x-y) F_{j}(s,y) dy ds\Big|^{p} dzdxdt \\
&\le \int_{0}^{T} \int_{{\mathbb R}^d} \int_{0}^{1} t^{p-1}\int_{0}^{t} \Big|{\sum_{j=1}^\infty} r_{j}(z) \int_{{\mathbb R}^d} P_{j}^{\alpha}(t-s,x-y) F_{j}(s,y) dy\Big|^{p} ds dzdxdt \\
&= \int_{0}^{T} t^{p-1}\int_{0}^{t} \int_{0}^{1} \int_{{\mathbb R}^d} \Big|{\sum_{j=1}^\infty} r_{j}(z) \int_{{\mathbb R}^d} P_{j}^{\alpha}(t-s,x-y) F_{j}(s,y) dy\Big|^{p} dx dzdsdt \\
&\le \int_{0}^{T} t^{p-1}\int_{0}^{t} \left({\sum_{j=1}^\infty} \left(\int_{{\mathbb R}^d} \Big|\int_{{\mathbb R}^d} P_{j}^{\alpha}(t-s,x-y) F_{j}(s,y) dy\Big|^{p} dx\right)^{1/p}\right)^{p} dsdt.
\end{split}
\end{equation}

\textbf{Step 2:} 
In order to estimate the double integral inside the sum, we shall use Young's inequality.
We claim that there exists a positive constant $c$ such that 
\begin{equation}
\label{s3:e3}
\int_{{\mathbb R}^d} |P_{j}^{\alpha}(t,x)| dx \ll 2^{j\alpha/p} \exp(-ct2^{j\alpha}).
\end{equation}
To see this, we use a change of variables to write 
\begin{align*}
P_{j}^{\alpha}(t,x) 
&= 2^{j\alpha/p} \int_{{\mathbb R}^d} e^{2\pi ix\cdot\xi} \widehat{\Psi}(2^{-j}\xi) e^{-t|\xi|^\alpha} d\xi \\
&= 2^{j\alpha/p} 2^{jd} \int_{{\mathbb R}^d} e^{2\pi i2^{j}x\cdot\xi} \widehat{\Psi}(\xi) e^{-t2^{j\alpha}|\xi|^\alpha} d\xi.
\end{align*}
From the observation 
\[(I-\Delta_\xi)e^{2\pi i2^{j}x\cdot\xi}=(1+4\pi^2 |2^{j}x|^2)e^{2\pi i2^{j}x\cdot\xi},\]
we can carry out repeated integration by parts so that we obtain for some positive number $c$ 
\[|P_{j}^{\alpha}(t,x)| \ll \frac{2^{j\alpha/p}2^{jd}\exp(-ct2^{j\alpha})}{(1+4\pi^2 |2^{j}x|^2)^{d+1}}.\]
Integrating with respect to $x$ gives the estimate \eqref{s3:e3}.
By Young's inequality, the last integral in \eqref{s3:e2} dominated by 
\[\int_{0}^{T} \int_{0}^{t} \left({\sum_{j=1}^\infty} 2^{j\alpha/p} \exp(-c(t-s)2^{j\alpha}) \Big(\int_{{\mathbb R}^d} |F_{j}(s,y)|^{p} dy\Big)^{1/p} \right)^{p} dsdt.\]

\textbf{Step 3:} 
Finally, we will estimate the following two integrals
\begin{align*}
I &:= \int_{0}^{T} \int_{0}^{t} \left(\sum_{(t-s)2^{j\alpha}\le1} 2^{j\alpha/p} \exp(-c(t-s)2^{j\alpha}) \Big(\int_{{\mathbb R}^d} |F_{j}(s,y)|^{p} dy\Big)^{1/p} \right)^{p} dsdt \\
II &:= \int_{0}^{T} \int_{0}^{t} \left(\sum_{(t-s)2^{j\alpha}>1} 2^{j\alpha/p} \exp(-c(t-s)2^{j\alpha}) \Big(\int_{{\mathbb R}^d} |F_{j}(s,y)|^{p} dy\Big)^{1/p} \right)^{p} dsdt.
\end{align*}
H\"older's inequality gives
\begin{align*}
&\left(\sum_{(t-s)2^{j\alpha}\le1} 2^{j\alpha/p} \exp(-c(t-s)2^{j\alpha}) \Big(\int_{{\mathbb R}^d} |F_{j}(s,y)|^{p} dy\Big)^{1/p} \right)^{p} \\
&\le \left(\sum_{(t-s)2^{j\alpha}\le1} 2^{j\alpha/(2p)} \Big(2^{j\alpha/2} \int_{{\mathbb R}^d} |F_{j}(s,y)|^{p} dy\Big)^{1/p} \right)^{p} \\
&\le \left(\sum_{(t-s)2^{j\alpha}\le1} 2^{j\alpha/(2(p-1))}\right)^{p-1} \sum_{(t-s)2^{j\alpha}\le1} 2^{j\alpha/2} \int_{{\mathbb R}^d} |F_{j}(s,y)|^{p} dy.
\end{align*}
Summing a geometric series gives
\[\left(\sum_{(t-s)2^{j\alpha}\le1} 2^{j\alpha/(2(p-1))}\right)^{p-1} \ll (t-s)^{-1/2}.\]
We change the order of integration and summation to obtain that  
\begin{align*}
I
&\le \int_{0}^{T} \int_{0}^{t} (t-s)^{-1/2} \sum_{(t-s)2^{j\alpha}\le1} 2^{j\alpha/2} \int_{{\mathbb R}^d} |F_{j}(s,y)|^{p} dy dsdt\\
&\le \int_{0}^{T} \left(\sum_{j=1}^{\infty} 2^{j\alpha/2} \int_{s}^{s+2^{-j\alpha}} (t-s)^{-1/2} dt \int_{{\mathbb R}^d} |F_{j}(s,y)|^{p} dy\right) ds\\
&\ll \int_{0}^{T} \left(\sum_{j=1}^{\infty} \int_{{\mathbb R}^d} |F_{j}(s,y)|^{p} dy\right) ds.
\end{align*}
For $q\le p$, the last quantity is dominated by 
\[\int_{0}^{T} \int_{{\mathbb R}^d} \left(\sum_{j=1}^{\infty} |F_{j}(s,y)|^{q}\right)^{p/q} dyds \le \norm{f}_{L^{p}F_{s}^{p,q}(Q)}^{p}.\]
Therefore, we get the desired result for the quantity $I$.

Similarly, we can estimate the quantity $II$.
H\"older's inequality gives
\begin{align*}
&\left(\sum_{(t-s)2^{j\alpha}>1} 2^{j\alpha/p} \exp(-c(t-s)2^{j\alpha}) \Big(\int_{{\mathbb R}^d} |F_{j}(s,y)|^{p} dy\Big)^{1/p} \right)^{p} \\
&= \left(\sum_{(t-s)2^{j\alpha}>1} 2^{j2\alpha/p} \exp(-c(t-s)2^{j\alpha}) \Big(2^{-j\alpha} \int_{{\mathbb R}^d} |F_{j}(s,y)|^{p} dy\Big)^{1/p} \right)^{p} \\
&\le \left(\sum_{(t-s)2^{j\alpha}>1} 2^{j2\alpha/(p-1)} \exp(-\widetilde{c}(t-s)2^{j\alpha})\right)^{p-1} \sum_{(t-s)2^{j\alpha}>1} 2^{-j\alpha} \int_{{\mathbb R}^d} |F_{j}(s,y)|^{p} dy,
\end{align*}
where the positive constant $\widetilde{c}$ depends only on $p$.
Summing a geometric series gives
\[\left(\sum_{(t-s)2^{j\alpha}>1} 2^{j2\alpha/(p-1)} \exp(-\widetilde{c}(t-s)2^{j\alpha})\right)^{p-1} \ll (t-s)^{-2}.\]
We change the order of integration and summation to obtain that  
\begin{align*}
II
&\le \int_{0}^{T} \int_{0}^{t} (t-s)^{-2} \sum_{(t-s)2^{j\alpha}>1} 2^{-j\alpha} \int_{{\mathbb R}^d} |F_{j}(s,y)|^{p} dy dsdt\\
&\le \int_{0}^{T} \left(\sum_{j=1}^{\infty} 2^{-j\alpha} \int_{s+2^{-j\alpha}}^{\infty} (t-s)^{-2} dt \int_{{\mathbb R}^d} |F_{j}(s,y)|^{p} dy\right) ds\\
&\ll \int_{0}^{T} \left(\sum_{j=1}^{\infty} \int_{{\mathbb R}^d} |F_{j}(s,y)|^{p} dy\right) ds.
\end{align*}
For $q\le p$, the last quantity is dominated by 
\[\int_{0}^{T} \int_{{\mathbb R}^d} \left(\sum_{j=1}^{\infty} |F_{j}(s,y)|^{q}\right)^{p/q} dyds \le \norm{f}_{L^{p}F_{s}^{p,q}(Q)}^{p}.\]
Therefore, we get the desired result for the quantity $II$.
This completes the proof of the theorem.

\end{document}